\documentclass[11pt]{amsart}
\usepackage{epsfig,amssymb, amsthm, graphicx, verbatim}
\newtheorem{theorem}{Theorem}[section]
\newtheorem{lemma}[theorem]{Lemma}
\newtheorem{proposition}[theorem]{Proposition}
\newtheorem{corollary}[theorem]{Corollary}
\newtheorem{conjecture}[theorem]{Conjecture}

\theoremstyle{definition}
\newtheorem{definition}[theorem]{Definition}

\newtheorem{remark}[theorem]{Remark}

\theoremstyle{remark}

\numberwithin{equation}{section}

\newcommand{\mcg}{\Gamma _{g,r}^k}
\newcommand{\bfz}{\mathbb{Z}}

\newcommand{\Mgnr}{{\Gamma }_{g,r}^k}
\newcommand{\Mgrn}{{\Gamma }_{g,r}^k}

\newcommand{\Mg}{{\Gamma }_{g}}

\newcommand{\bfc}{{\mathbb {C}}}
\newcommand{\bfr}{{\mathbb {R}}}

\newcommand{\bdy}{\partial}

\begin{document}

\newenvironment{prooff}{\medskip \par \noindent {\it Proof}\ }{\hfill
$\square$ \medskip \par}
    \def\sqr#1#2{{\vcenter{\hrule height.#2pt
        \hbox{\vrule width.#2pt height#1pt \kern#1pt
            \vrule width.#2pt}\hrule height.#2pt}}}
    \def\square{\mathchoice\sqr67\sqr67\sqr{2.1}6\sqr{1.5}6}
\def\pf#1{\medskip \par \noindent {\it #1.}\ }
\def\endpf{\hfill $\square$ \medskip \par}
\def\demo#1{\medskip \par \noindent {\it #1.}\ }
\def\enddemo{\medskip \par}
\def\qed{~\hfill$\square$}

 \title[Contact 3-manifolds with infinitely many Stein fillings]
 {Contact 3-manifolds with infinitely many Stein fillings}

 \author{Burak Ozbagci}
 \author{Andr\'{a}s I. Stipsicz}

 \address{Department of Mathematics, Michigan State University, 
 East Lansing, MI, 48824; current address: Dept. of Math., Koc University,
Istanbul, Turkey} 
 \email{bozbagci@math.msu.edu}

 \address{A. R\'enyi Institute of Mathematics, Hungarian Academy of Sciences, 
 Budapest, Hungary and Department of Mathematics, Princeton University, 
 Princeton, NJ, 08544} 
 \email{stipsicz@math-inst.hu and stipsicz@math.princeton.edu}
 \date{\today}

\begin{abstract}
Infinitely many contact 3-manifolds each admitting infinitely many,
pairwise non-diffeomorphic Stein fillings are constructed. We use Lefschetz
fibrations in our constructions and compute their first homologies to 
distinguish the fillings.
\end{abstract}
\maketitle

\section{Introduction}\label{first}

A complex surface $V$ is {\em Stein\/} if it admits a proper
holomorphic embedding $f\colon V\to \bfc ^n$ for some $n$. For a
generic point $p\in \bfc ^n$ consider the map $\varphi \colon V\to
\bfr $ defined by $\varphi (z)=\vert \vert z-p \vert \vert ^2$. For a
regular value $a\in \bfr$ the level set $M=\varphi ^{-1}(a)$ is a
smooth 3-manifold (oriented as $\partial \varphi ^{-1}([0,a])$) with a
distinguished 2-plane field $\xi = TM\cap iTM \subset TV$. It turns
out that $\xi $ defines a {\em contact structure \/} on $M$ (for more
about contact structures see \cite{Ae, E}) and $S= \varphi
^{-1}([0,a])$ is called a {\em Stein filling\/} of $(M, \xi
)$. Topological properties of Stein fillings (and slightly more
generally, of strong symplectic fillings) are in the focus of current
research. Based on work of Eliashberg \cite{Eli}, McDuff \cite{M} and
Lisca \cite{L} we know, for example, that for the lens space $L(p,q)$
equipped with a specific contact structure $\xi _{(p,q)}$ there is a
finite list $\{ S_{p,q}(n)\mid n=1, \ldots , n_{p,q}\}$ of 4-manifolds
such that any Stein filling of $(L(p,q), \xi _{(p,q)})$ is
diffeomorphic to some $S_{p,q}(n)$. Similar finiteness results have
been verified for simple and simple elliptic singularities \cite{OO1, OO2},
and for homeomorphism types of Stein fillings of the 3-torus $T^3$
\cite{S}.  Based on these examples it was anticipated that a contact
3-manifold $(M, \xi )$ admits only finitely many non-diffeomorphic
Stein fillings.  Genus-3 fibrations found by I. Smith \cite{sm2}
indicated that such an expectation is too ambitious in general.

In the following we show the existence of
an infinite family of contact 3-manifolds each
admitting infinitely many non-diffeomorphic Stein fillings.

\begin{theorem}\label{main}
For all $g\geq 2$ there is a contact 3-manifold $(M_g, \xi _g)$
such that $M_{g}$ is diffeomorphic to $M_{g'}$ iff $g=g'$ and 
each $(M_g, \xi _g)$ admits infinitely many pairwise non-diffeomorphic 
Stein fillings.
\end{theorem}  

In fact, $M_g$ can be given as the boundary of the plumbing of the
disk bundle 
over a genus-$g$
surface with Euler number zero and the  disk bundle 
over a sphere with Euler number $2$. 
In the proof of Theorem~\ref{main} we use Lefschetz fibrations 
to construct the Stein fillings and we compute their first homology
groups to distinguish them.
It would be most desirable to find infinitely many distinct fillings
with trivial fundamental group; we hope to return to this point later. 
Regarding finiteness of Stein fillings, we conjecture the following:
\begin{conjecture}\label{finite}
Let $\chi (X)$ denote the Euler characteristic of the
compact manifold $X$.
Then the set ${\mathcal {C}}_{(M, \xi )}=\{ \chi (S) \mid S
{\mbox{ is a Stein filling of }} (M, \xi )\}$ is finite.
\end{conjecture}

\noindent
{\bf {Acknowledgement}}: The authors would like to thank Selman Akbulut,
Mustafa Korkmaz and
John Etnyre for many helpful and inspiring discussions. The second author 
was partially supported by OTKA.

\section{Fiber sums of Lefschetz fibrations}\label{second}

Let $\Sigma$ be a compact, oriented and connected surface of genus $g$ 
with $r$ marked points and $k$ boundary components. The 
{\em mapping class group} $\mcg$ of 
$\Sigma$ consists of the isotopy classes of 
orientation-preserving self-diffeomorphisms of $\Sigma$, which are 
identity on each boundary component and preserve the set of marked points. 
(The groups $\Gamma_{g,r}^0 $, $\Gamma_{g,0}^k$ and 
$\Gamma_{g,0}^0$ will be abbreviated by  $\Gamma_{g,r}$, $\Gamma_{g}^k$ and 
$\Mg$, respectively.) 
We say that a $\Sigma$-bundle $P$ over $S^1$
has monodromy $h\in  \Mgnr $ iff $P$ is diffeomorphic to
$ (\Sigma \times I ) / \;(h(x),0) \sim (x,1)$.

\begin{definition}\label{lefi}
Let $X$ be a compact, connected, oriented, smooth four-manifold. 
A Lefschetz fibration on $X$ is a smooth 
map $\pi \colon X \to B$, where $B$ is a 
compact, connected, oriented surface,
$ {\pi}^{-1} (\partial B) = \partial X$, furthermore 
each critical point of $\pi$ lies in int$X$ and has an 
orientation-preserving local coordinate chart on 
which $\pi ( z_1 , z_2 ) = z_1^2 + z_2^2$.    
\end{definition}

It follows that  $\pi$ has only finitely many critical points 
and  removing the 
corresponding singular fibers turns a Lefschetz fibration into a 
fiber bundle with a connected base space. Consequently 
all but finitely many 
fibers of a Lefschetz fibration are smooth, 
compact and oriented 
surfaces, all of which having the same diffeomorphism type
of a closed genus-$g$ surface for some $g$.
We will assume that there is at most one critical point 
on each fiber. A Lefschetz 
fibration is called
{\em relatively minimal} if no fiber contains
an embedded 2-sphere of self-intersection number $-1$.
Each critical point of a Lefschetz fibration corresponds to 
an embedded circle in a nearby regular fiber
called a {\em vanishing  cycle}, and the singular 
fiber is obtained by collapsing the vanishing cycle to a point. 
The boundary of a regular neighborhood of a singular fiber 
is a surface bundle over the circle. In fact, a singular fiber 
can be described by the monodromy of this surface bundle which turns out to be
a right-handed Dehn twist along the corresponding vanishing cycle. 
Once we fix an identification of $\Sigma$ with the fiber over 
a base point of $B$, the topology of the Lefschetz fibration is 
determined by its {\em monodromy representation} 
$\Psi \colon {\pi}_1 (B- \{ {\mbox {critical values}} \}) \to \Gamma _g .$
In case $B=D^2$ the monodromy along $\partial D^2=S^1$ is called the {\em total
monodromy\/} of the fibration; according to the above said it is the 
product of right-handed Dehn twists corresponding to the singular fibers.
A Lefschetz fibration over 
$S^2$ can be decomposed into two Lefschetz fibrations over $D^2$, 
one of which is trivial; consequently a Lefschetz fibration over 
$S^2$ is determined by a relator in the mapping class group. Conversely, 
given a product of right-handed  Dehn twists in the mapping class group 
we can construct the corresponding Lefschetz fibration over $D^2$, and 
if the given product of right-handed  Dehn twists 
is isotopic to identity (and $g\geq 2$) then the fibration 
extends uniquely over $S^2$.
The monodromy representation also provides a handlebody decomposition
of a Lefschetz fibration over $D^2$: we attach 2-handles to $\Sigma \times D^2$
along the vanishing cycles with framing $-1$ relative to the framing the
circle inherits from the fiber.
(For a more detailed introduction to the theory of Lefschetz fibration see
\cite{GS}.)
 
Let $X \to S^2$ 
be a Lefschetz fibration with generic fiber $\Sigma$ and 
let $\gamma_1, \gamma_2, \cdots , \gamma_s$ denote the 
vanishing cycles of this fibration. 
Assume that $X \to S^2$ 
admits a section, i.e., there is $\sigma \colon S^2 \to X$ with
$\pi \circ \sigma = {\mbox {id}}_{S^2}$. 
The following results are standard.

\begin{lemma}\label{fundamental} 
The first homology group $H_1 (X; \bfz)$ is the quotient of 
$H_1 ( \Sigma ; \bfz )$ by the subgroup generated by the 
homology classes of the vanishing cycles.\qed 
\end{lemma}

\begin{lemma}\label{fibersum}
Let $X \sharp_{f} X $ denote the fiber sum of X with 
itself by a self-diffeomorphism $f$ of the generic fiber $\Sigma$.
Then $X \sharp_{f} X $ is a Lefschetz fibration with vanishing cycles 
$\gamma_1, \gamma_2, \cdots , \gamma_s, f({\gamma}_1), 
f({\gamma}_2), \cdots , f({\gamma}_s) $.  \qed
\end{lemma}

\begin{remark}\label{conjugate}
If $\alpha$ is a simple closed curve on $\Sigma$,
$t_{\alpha }$ is the corresponding Dehn twist and $f$ is an 
orientation-preserving self-diffeomorphism of $\Sigma$, then 
$f t_{\alpha} f^{-1} = t_{f(\alpha)}$ in $\Mgrn$.
\end{remark}

\begin{figure}[ht]
  \begin{center}
     \includegraphics{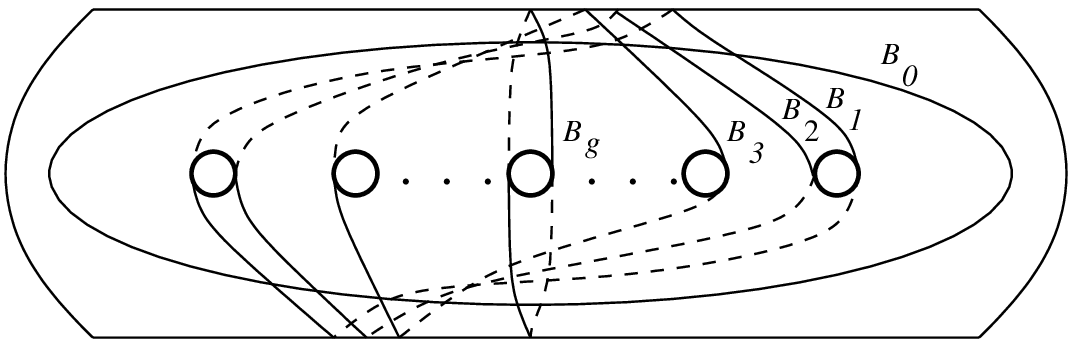}
   \caption{}\label{word}
    \end{center}
  \end{figure}

\begin{figure}[ht]
  \begin{center}
     \includegraphics{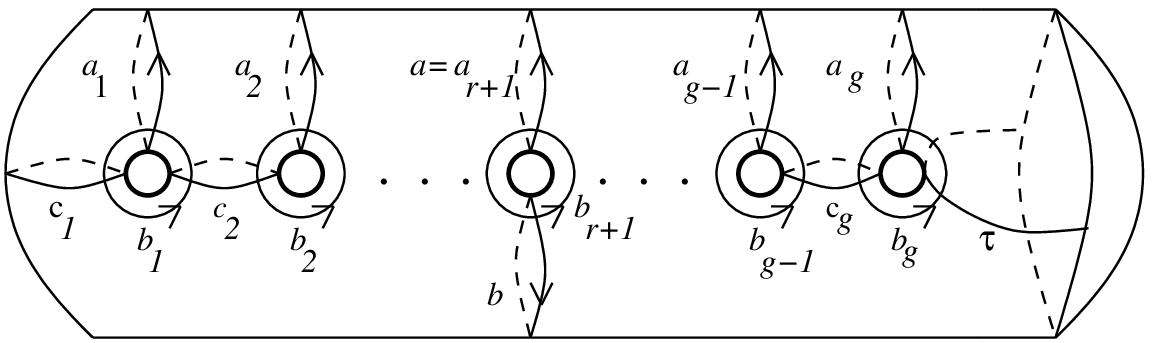}
   \caption{}\label{s}
    \end{center}
  \end{figure}

After this short introduction we begin our construction with a
description of a set of words in the mapping class groups which was
discovered by Korkmaz \cite{ko}.  We focus on the odd genus case, for
$g$ even see Remark~\ref{eveng}.  For $g = 2r +1 \geq 3$ 
the following relation holds in $\Mg$:
$$W_g = (t_{B_0} t_{B_1} t_{B_2} \cdots t_{B_g} t_a^2 t_b^2 )^2 = 1  $$
where  $B_0 , B_1 , \cdots , B_g$ are shown in 
Figure~\ref{word} and  $a$ and $b$ in 
Figure~\ref{s}. 
Let $X_g \to S^2 $ denote the Lefschetz fibration 
corresponding to the relator $W_g $ in $\Mg$. Now consider the fiber sum 
of $X_g$ with itself using 
the diffeomorphism $ t_{a_1}^n $ and denote the result by $X_g(n)$. 
(The simple closed curve $a_1$ is depicted in Figure~\ref{s}.)
Then by Lemma~\ref{fibersum}
the 4-manifold $X_g(n)$ comes with a Lefschetz fibration 
$X_g(n) \to S^2$ of global monodromy 
$W_g(n) = W_g W_{g}^{t_{a_1}^n} = 1 .$ 
(Here $W_{g}^{t_{a_1}^n}$ means the conjugate of $W_{g}$ by $t_{a_1}^n$.) 
Notice that by Remark~\ref{conjugate} $W_{g}^{t_{a_1}^n}$ and hence 
$W_g(n)$ are products of right-handed  Dehn twists.

It is routine to check that the curves 
$ c_1, b_1 , c_2 , b_2 , \cdots , c_g , b_g$ and $a_2$ 
(depicted in Figure~\ref{s}) are fixed by the 
product $(t_{B_0} t_{B_1} t_{B_2} \cdots t_{B_g} t_a^2 t_b^2 )^2$ 
on the surface with one puncture. This shows that 
$(t_{B_0} t_{B_1} t_{B_2} \cdots t_{B_g} t_a^2 t_b^2 )^2 = 1 $ holds in 
$\Gamma_{g,1}$.  Next we determine the
element $W_g$ considered in $\Gamma _g ^1$. It is known that
$\ker \varphi =\{ \Delta _g ^n\} \cong \bfz $ for 
the natural homomorphism $\varphi \colon 
\Gamma _g ^1 \to \Gamma _{g,1}$ collapsing a boundary 
circle to a marked point.  
(Here ${\Delta}_g$ denotes the Dehn twist 
along a curve parallel to the boundary circle.) Hence $W_g=\Delta _g ^n$
for some $n \in \bfz$ follows from the above discussion.

\begin{lemma}\label{boundary} 
$ (t_{B_0} t_{B_1} t_{B_2} 
\cdots t_{B_g} t_a^2 t_b^2 )^2 = {\Delta}_g $ holds 
in ${\Gamma}_g^{1} $.
\end{lemma}

\begin{proof}
We depicted $t_b^2 t_{B_0} t_{B_1} t_{B_2} 
\cdots t_{B_g} t_a^2 t_b^2 (\tau) $ and $t_{B_{g-1}}^{-1} t_{B_{g-2}}^{-1}  
\cdots t_{B_0}^{-1} {\Delta}_g (\tau) $ in Figures~\ref{bound1} 
and \ref{bound2}, respectively, where $\tau$ is shown 
in Figure~\ref{s}. One can use induction to obtain these figures. 
For example, the $W$-shaped part in the middle of Figure~\ref{bound2} 
is obtained by moving the legs of $W$ from outside to the 
inside of the holes when we apply 
$ t_{B_{i}}^{-1}  t_{B_{i-1}}^{-1}$. Now it is easy to see 
that $$ t_a^2 t_b^2 t_{B_0} t_{B_1} t_{B_2} 
\cdots t_{B_g} t_a^2 t_b^2 (\tau) =  t_{B_g}^{-1} t_{B_{g-1}}^{-1}  
\cdots t_{B_0}^{-1} {\Delta}_g (\tau) $$ which proves the lemma. 
\begin{figure}[ht]
  \begin{center}
     \includegraphics{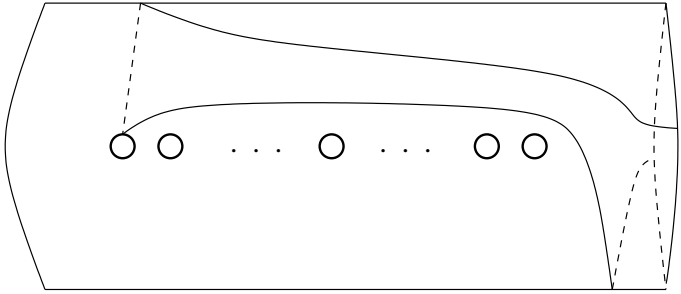}
   \caption{ $t_b^2 t_{B_0} t_{B_1} t_{B_2} 
\cdots t_{B_g} t_a^2 t_b^2 (\tau) $}\label{bound1}
    \end{center}
  \end{figure}
\begin{figure}[ht]
  \begin{center}
     \includegraphics{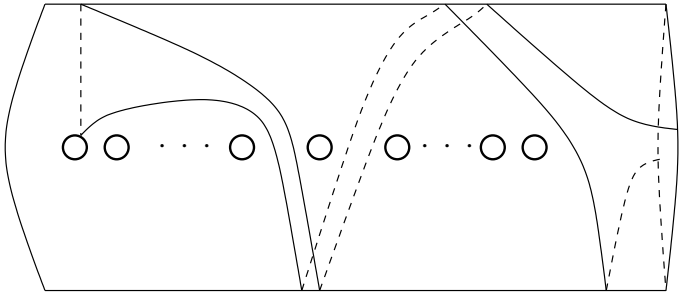}
   \caption{$ t_{B_{g-1}}^{-1}  t_{B_{g-2}}^{-1}
\cdots t_{B_0}^{-1} {\Delta}_g (\tau) $}\label{bound2}
    \end{center}
  \end{figure}
\end{proof}

\begin{corollary}\label{section} 
The Lefschetz fibration $X_g(n) \to S^2 $ admits a section 
$\sigma_g(n)$ with 
self-intersection number $-2$.  
\end{corollary}

\begin{proof}
We saw that
the relator $W_g$ admits a lift under the obvious map 
$\Gamma_{g,1} \to \Gamma_{g}$ defined by forgetting the 
marked point. Once such a lift is fixed, it gives a 
section $\sigma_g$ of the Lefschetz fibration $X_g \to S^2 $.
Similarly we get a section $\sigma_g(n)$ of $X_g(n) \to S^2 $. 
Lemma~\ref{boundary} proved that 
${\Delta}_g = W_g $ and thus ${\Delta}_g^2 = W_g(n)$ in  
$\Gamma_g^1$ since $t_{a_1}^{n} {\Delta}_g  t_{a_1}^{-n} = {\Delta}_g$.  
It is known (see \cite{sm1}) that if for the relator $w$ we have
$w=\Delta _g ^k$ in $\Gamma _g ^1$ then the Lefschetz fibration 
$X\to S^2$ given by $w$ admits a section of square $-k$.
This observation concludes the proof. 
\end{proof}

Since our Lefschetz fibrations admit sections,
Lemmma~\ref{fundamental} 
can be  applied to compute their first homology groups:
\begin{lemma}
The first homology group $H_1(X_g(n); \bfz )$ of $X_g(n)$ is isomorphic 
to  ${\bfz}^{g-2} \oplus  {\bfz}_n$. 
\end{lemma}

\begin{proof}
In the following we will denote the homology classes of curves 
by the same letters as we denote the curves. Let $a_1 , b_1 , 
a_2 , b_2, \cdots , a_g , b_g$ denote the standard generators of 
the first homology group of the fiber, as depicted in Figure~\ref{s}. 
It is easy to see that $B_i^{t^n_{a_1}}=B_i$ for $i\geq 2$,
and for an appropriate choice of orientation on $B_i$, in homology we have 
$$ B_{g} = a + b_{r+1} + b, \mbox{\ \ \ \ \ \ }  B_0 = b_1 + b_2 + \cdots + b_g , $$
$$B_i^{t^n_{a_1}}=B_i+na_1 \mbox{\ \ \ for\ \ } i=0,1,$$
$$ B_{2i-1} = b_{i} + B_{2i} + b_{g-(i-1)} \mbox{\ \ \ and} $$ 
$$ B_{2i} = a_i - a_{i+1} + B_{2i+1} + a_{g-(i-1)} - a_{g-i} $$
for $i = 1, \cdots ,r$, where $g= 2r+1$.
Thus by Lemmas~\ref{fundamental} and \ref{fibersum} the homology
$H_1 ( X_g(n); \bfz )$ is an 
(abelian) group generated by $ a_1, b_1 , a_2 , b_2 , \cdots , 
a_g , b_g $ with relations 
$$ a_{r+1} = b_{r+1} = 0, \ \ \ 
a_i + a_{g-(i-1)} = 0, \ \ \ b_i + b_{g-(i-1)} = 0 {\mbox { and }} n a_1 =0.$$ 
Now Lemma~\ref{fundamental} implies the result.
\end{proof}

In fact, for any given integer $ 0 \leq  k \leq g-2 $ and $n \geq 1$, 
we can construct a Lefschetz fibration over $S^2$ with first homology 
group isomorphic to ${\bfz}^k \oplus {\bfz}_n$, by taking 
twisted fiber sums of more copies of $X_g$. We could use 
any of these fibrations in the rest of the paper as well.

\begin{figure}[ht]
  \begin{center}
     \includegraphics{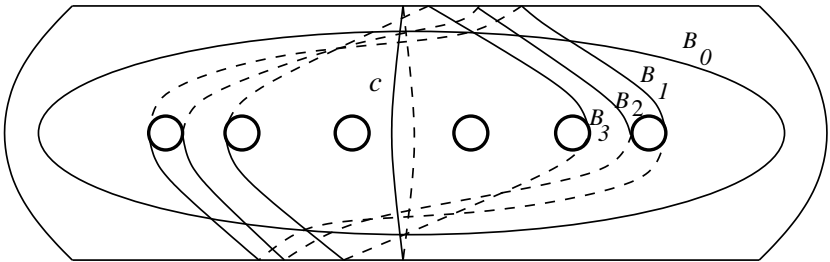}
   \caption{Vanishing cycles in the even $g$ case}\label{even}
    \end{center}
  \end{figure}

\begin{remark}\label{eveng}
If $g=2r\geq 2$ is even, similar results can be obtained for
fibrations defined by the relator $W_g= (t_{B_0} t_{B_1} t_{B_2}
\cdots t_{B_g} t_c)^2$ (see Figure~\ref{even}).  In the exact same
manner we can construct $X_g$ and the twisted fiber sums $X_g(n)\to
S^2$ with sections of square $-2$. The homology computation applies
without essential change and provides $H_1(X_g (n); \bfz ) =\bfz
^{g-2}\oplus \bfz _n$. (The observation preceding this remark also has
its natural extension to the even $g$ case.)  The only notable
difference is that for $g$ even the relator $W_g$ contains
homologically trivial vanishing cycles as well --- cf. the curve
denoted by $c$ on Figure~\ref{even}.
\end{remark}

\section{Stein structures on Lefschetz fibrations with boundary}\label{third}

\begin{definition}  
Let $X$ be a compact 4-manifold with boundary. The map 
$\pi \colon X \to B$ is a Lefschetz fibration if $\pi $ satisfies 
the requirements of Definition~\ref{lefi} with the possible exemption
of allowing $\pi ^{-1}(\partial B)\neq \partial X$, and assuming that
$\pi $ is a fiber bundle map away from the singular fibers.
(Notice that this definition allows the fiber $\pi ^{-1}(t)$ to have boundary.)
In the following, LFB (Lefschetz fibration with bounded fibers)
will denote a relatively minimal Lefschetz fibration over 
$D^2$ whose generic fiber is a surface with nonempty boundary. 
\end{definition}

\begin{theorem}[\cite{ao1, eh, lp}] \label{palf}
Let $f\colon X \to S^2 $ be a Lefschetz fibration with a section $\sigma$ 
and let $\Sigma$ denote 
a regular fiber of this fibration. 
Then $S = X - {\mbox {int }}\nu ( \Sigma \cup \sigma ) $ is a Stein filling 
of its boundary equipped with the induced (tight) contact structure, 
where $ \nu ( \Sigma \cup \sigma )$ denotes a regular 
neighborhood  of $\Sigma \cup \sigma$ in $X$.
\end{theorem}
We only sketch the proof of the above theorem, since 
it was already proved in \cite{lp, ao1} (with the extra assumption 
of having only homologically essential vanishing cycles in the fibration),
and also appeared in \cite{eh}.

\begin{prooff}{\em of Theorem~\ref{palf} (sketch).}
When we remove  ${\mbox {int }}\nu ( \Sigma \cup \sigma )$ 
from $X$, it is clear 
that we get a LFB. Thus 
$S = X - {\mbox {int }}\nu ( \Sigma \cup \sigma ) $ admits a standard 
handlebody decomposition: It is obtained by a
sequence of steps of
attaching 2-handles
$S_{0}= D^2 \times F  \leadsto  S_{1} \leadsto S_{2} \cdots \leadsto
S_{n}=S \;$, where each $S_{i-1}$ is a LFB
and $S_{i}$ is obtained from  $S_{i-1}$ by attaching a
2-handle along a simple closed curve
$\gamma$ lying on a fiber ${\Sigma}^{\prime} 
\subset \partial S_{i-1}$. Furthermore this handle is
attached to $\gamma$ with the framing $k-1$, where $k$ is the
framing induced from the
surface 
${\Sigma}^{\prime} $ (see text after Definition~\ref{lefi}). 
Inductively we assume that $ S_{i-1} $ has a Stein
structure, with a
convex fiber $ {\Sigma}^{\prime} \subset \partial S_{i-1}$. By
the ``Legendrian realization principle" of
\cite{ho} (pp 323-325)
applied to the
homologically essential vanishing cycle $\gamma$, 
after an isotopy of $({\Sigma}^{\prime} 
, \gamma )$,  $\gamma$ becomes a Legendrian curve and $k$
can be taken to be its
Thurston-Bennequin framing. 
If $\gamma $ is homologically inessential, the above principle does not
apply verbatim. In this case choose a homologically essential simple closed
curve in the component of the complement of $\gamma$  which is disjoint from
$\partial \Sigma '$. Isotope it to Legendrian position (by the above argument)
 and using the local model of the contact structure along this Legendrian
knot $\delta $ isotope $\Sigma '$ to introduce new dividing curves
parallel to $\delta $. Now the Legendrian realization principle applies 
to $\gamma $ since after this modification (called a {\em fold\/})
each component of its complement intersects the set
of dividing curves nontrivially. 
(For the relevant definitions see \cite{E, ho}.)
Now a celebrated result of Eliashberg \cite{el, GS} provides an extension
of the Stein structure of $S_{i-1}$ to $S_i$.
\end{prooff}

\begin{remark}
Notice that if $g$ is odd then all vanishing cycles are homologically 
essential, hence in this case no fold is necessary. In the even $g$ case
our examples contain homologically inessential vanishing cycles ($c$ on 
Figure~\ref{even}), hence application of folds is necessary.
\end{remark}

\begin{definition} 
Let $M$ be a closed $3$-manifold. An {\em open book\/} decomposition of 
$M$ is a pair $(K,\pi )$ consisting of a (fibered) link $K$,
 called the {\em binding\/}, and a 
fibration $\pi \colon  M-K \to S^1$ such that each fiber is a Seifert 
surface for $K$. 
(The fibers are also called the {\em pages\/} of the open book.)
Let ${\Sigma}^{\prime}$ denote a surface with boundary. 
The closed 3-manifold 
$$M = (({\Sigma}^{\prime} \times I ) / \;(h(x),0) \sim (x,1))   
\cup_{\partial} (\bdy {\Sigma}^{\prime} \times D^2)$$ is 
canonically decomposed as  
an open book with binding $\bdy {\Sigma}^{\prime} $, 
page ${\Sigma}^{\prime} $ and 
monodromy $h$. 
\end{definition}

It is easy to see that the boundary of a LFB has a canonical
open book decomposition induced from the fibration. 
The pages are the fibers in the boundary of the LFB, while the monodromy
of the open book is just the total monodromy of the LFB along 
$S^1=\partial D^2$. The binding is simply the boundary of the
central fiber $\pi ^{-1}(0)$.

We now turn our attention to contact structures on $3$-manifolds. 
In \cite{tw}, Thurston and Winkelnkemper constructed 
a contact structure associated to a given 
open book decomposition of a closed, orientable $3$-manifold. 
Recently Giroux \cite{Gi} 
refined this construction by showing that an open book $(K, \pi )$
supports a unique (up to isotopy) {\em compatible\/} contact structure
$\xi$, where a contact structure $\xi $ is called
compatible with an open book $(K, \pi )$ if 
there is a contact 1-form $\alpha$ for $\xi $ such that $d\alpha $ 
is a volume form on each fiber of $\pi $ and  
the binding $K$ is a transverse link in $(M,\xi )$ (oriented as
the boundary of a page). In fact, Giroux established  a bijection
between the set of isotopy classes of contact structures on a closed
3-manifold $M$ and open books supported by $M$; here two open books are
considered to be equivalent if they can be joined by a sequence of
positive stabilization/destabilization.  
In addition, Giroux proved that $\xi $
is Stein fillable if the monodromy 
of the compatible open book decomposition 
can be expressed as a product of right-handed  Dehn twists.
Notice that the boundary of a LFB admits
two contact structures --- one induced by the Stein structure of the LFB and
the other given by the open book through the total monodromy of the
Lefschetz fibration. 
The above  result of Giroux \cite{Gi} essentially states that these two 
structures are isotopic.
(A slightly different proof of the same statement --- resting on 
results of Torisu~\cite{To} --- has been given by Gay~\cite{Gay}.)

\begin{proposition}[\cite{eh, Gay, Gi}]\label{contact} 
The contact structure 
on the boundary of a LFB 
(induced from the Stein structure) is isotopic to the 
contact structure associated to the boundary of this 
LFB considered as an open book. In particular, 
any two LFB's 
bounding the same positive 
open book decomposition are Stein 
fillings of the same 
contact structure. \qed
\end{proposition}

\begin{remark}
In proving this proposition one first checks it for the trivial
LFB with no critical points (hence with total monodromy equal to id).
In this case the Lefschetz fibration is built using 1-handles only.
Then an induction on the number of 2-handles together with the 
observation that the handle attachment 
(which is surgery on the 3-manifold level) yields
 an open book compatible with the
new contact structure proves the result. (This last observation
only involves a local check along a page of the open book.)
\end{remark}

In the following we apply the above discussion to the 
Lefschetz fibrations $X_g(n) \to S^2 $ we constructed in 
Section~\ref{second}. 
Let $\Sigma_g(n)$ denote a regular fiber of the 
Lefschetz fibration  $X_g(n) \to S^2 $ with section  $\sigma_g(n)$. 
Let $S_g(n)$ denote 
$X_g(n) - {\mbox {int }}\nu ( \Sigma_g(n) \cup \sigma_g(n) )$. 

\begin{lemma}\label{homology}
The first homology group $H_1 ( S_g(n) ; \bfz )$ of $S_g(n)$ is 
isomorphic to $H_1(X_g(n); \bfz ) =  {\bfz}^{g-2} \oplus {\bfz}_n$. 
\end{lemma}

\begin{proof}
$H_1(X_g(n); \bfz ) \cong H_1(S_g(n); \bfz )$ 
follows from the observation that the normal circle to the fiber 
is homologous to a multiple of the normal circle to the section
(shown by a push-off of the section), which bounds a punctured 
fiber, hence is zero in $H_1(S_g(n); \bfz )$.
\end{proof}

\begin{prooff}{\em of Theorem~\ref{main}.}
Let $M_g(n)$ denote the boundary of $S_g(n)$ with the contact structure
$\xi _g (n)$ induced by the open book provided by $S_g(n)$ as a LFB.
Notice that $M_g(n)=-\partial \nu (\Sigma \cup \sigma )$, hence its
diffeomorphism type does not depend on $n$. 
Moreover, the total monodromy of $S_g(n)$ is equal to $\Delta _g^2\in \Gamma
_g ^1$, hence the open book, and therefore the contact structure $\xi _g (n)$
is independent of $n$. Let $(M_g, \xi _g)$ denote $(M_g(n), \xi _g (n))$.
According to Proposition~\ref{contact} the Stein structures on 
$S_g(n)$ ($n=1,2,\ldots $) all provide Stein fillings of $(M_g, \xi _g)$,
and using Lemma~\ref{homology} we conclude that these are pairwise
non-diffeomorphic fillings.

Since $H_1 (M_g ; \bfz) = {\bfz}^{2g}$, we also see that $M_g$ is
diffeomorphic to $M_{g^{\prime}} $ if and only if $g = g^{\prime}$,
therefore the proof of Theorem~\ref{main} is complete.  Since $\nu
(\Sigma \cup \sigma )$ is just the plumbing of the disk bundle over a
genus-$g$ surface with Euler number 0 and the disk bundle over a
sphere with Euler number $-2$, we get that $M_g$ (which is $\partial
\nu ( \Sigma \cup \sigma )$ with the opposite orientation) is the
boundary of a similar plumbing (now with Euler number $2$), as stated
in Section~\ref{first}.
\end{prooff}

\begin{remark}
The above construction just gave different factorizations of the 
mapping class $\Delta _g^2 \in \Gamma _g ^1$ into products of 
right-handed Dehn twists. In fact, any such factorization of an element 
$h\in \Gamma _g ^r$ ($r>0$) provides a Stein filling of the contact 
structure induced by $h$, and the topological argument (computing the
first homology) helped us to distinguish the various factorizations of
$\Delta _g ^2$. In the light of this observation we can weaken 
Conjecture~\ref{finite} to
\end{remark}

\begin{conjecture}\label{bound}
For $h\in \Gamma _g ^r$ $(r>0)$ there exists a constant $C_h$
such that if $h=t_1\cdot \cdot \cdot t_n$ with $t_i$ right-handed
Dehn twists then $n\leq C_h$.
\end{conjecture}

The following observation can serve as an evidence for Conjecture~\ref{bound}.

\begin{lemma} For $g$ odd the Euler 
characteristic $\chi ( S_g(n)) $ and the 
signature $\sigma (  S_g(n) )$ of $S_g(n)$ 
are given by $11$ and $ -16$, respectively.   
\end{lemma}

\begin{proof}
For $g$ odd the 
Stein filling $S_g(n)$ admits a LFB with $2g +10$ singular fibers. 
Thus $$ \chi ( S_g(n)) = \chi(D^2) \chi( {\mbox {fiber}}) + 
2g +10 = 2- 2g -1 + 2g +10 = 11 .$$ 
The signature $\sigma  (  X_g(n) ) = -16 $ since 
$\sigma (  X_g ) = -8 $ \cite{ko} and the signature 
is additive under fiber sum. It implies that  
$\sigma (  S_g(n) ) = - 16 $  since we remove a piece with 
zero signature to get $S_g(n)$ from $X_g(n)$.
\end{proof}
(For $g$ even a similar computation gives $\chi (S_g(n))=5$ 
and $\sigma (S_g(n))=-8$.)
Notice that the length of the factorization $h=t_1\cdot \cdot \cdot
t_n$ into right-handed Dehn twists is intimately related to the
Euler characteristic of the Stein filling induced by
$t_1\cdot \cdot \cdot t_n$.

\end{document}